\documentclass[12pt]{amsart}
\usepackage{amsmath,amssymb,amsbsy,amsfonts,amsthm,latexsym,amsopn,amstext,
                                                       amsxtra,euscript,amscd}
\begin{document}

\title[Riemann Hypothesis may be proved by induction]{Riemann Hypothesis \\may be proved by induction}

\author[R. M.~Abrarov]{R. M.~Abrarov}
\address{University of Toronto, Canada}
\email{rabrarov@physics.utoronto.ca}

\author[S. M.~Abrarov]{S. M.~Abrarov}
\address{York University, Toronto, Canada}
\email{abrarov@yorku.ca}

\date{\today}

\begin{abstract} 

The transformations of the sum identities for generalized harmonic and oscillatory numbers, obtained earlier in our recent report [1], enable us to derive the new identities expressed in terms of the corresponding square roots of $x$. At least one of these identities may be applied to prove the Riemann Hypothesis by induction. Additionally using this approach, the new series for Euler's constant $\gamma$ has been found.
\end{abstract}

\maketitle

Let us start discussion from the sum identity
\begin{equation} \label{Eq_1} 
\sum _{i=1}^{x}\frac{1}{i^{s} }  M_{\frac{x}{i} } \left(s\right)=1 
\end{equation} 
obtained in our recent work [1] for generalized oscillatory numbers in power \textit{s}
$$
M_{x} \left(s\right)=\sum _{k=1}^{x}\frac{\mu _{k} }{k^{s} }  ,
$$
where $\mu _{k} $ is the M\"obius function.
Rearranging $M_{\frac{x}{i} } \left(s\right)$ in terms of $\frac{\left(x^{2} /j\right)}{i} $, where $j\ge x$ is an integer variable, the equation \eqref{Eq_1} can be rewritten as
\begin{equation}\label{Eq_2}
\sum _{i=1}^{x}\frac{1}{i^{s} }  M_{\frac{\left(x^{2} /j\right)}{i} } \left(s\right)=1,\qquad j\ge x .
\end{equation}
Multiplying both parts of \eqref{Eq_1} by $\frac{1}{j^{s} } $ and taking the sum over index \textit{j} from $x+1$ up to $x^{2} $, the right hand side can be represented as a difference of the generalized harmonic numbers in power \textit{s} [2] for $x^{2} $ and $x$, i.e.:
\begin{equation}\label{Eq_3}
\sum _{j=x+1}^{x^{2} }\frac{1}{j^{s} } \sum _{i=1}^{x}\frac{1}{i^{s} }   M_{\frac{\left(x^{2} /j\right)}{i} } \left(s\right)=\sum _{j=x+1}^{x^{2} }\frac{1}{j^{s} }  =H_{x^{2} } \left(s\right)-H_{x} \left(s\right). 
\end{equation}
On the other hand we can rewrite \eqref{Eq_3} as
\begin{equation} \label{Eq_4} 
\begin{aligned} 
\sum _{i=1}^{x}\frac{1}{i^{s} }  \sum _{j=x+1}^{x^{2} }\frac{1}{j^{s} }  M_{\frac{\left(x^{2} /i\right)}{j} } \left(s\right)&=\sum _{i=1}^{x}\frac{1}{i^{s} }  \left(1-\sum _{j=1}^{x}\frac{1}{j^{s} }  M_{\frac{\left(x^{2} /i\right)}{j} } \left(s\right)\right) \\ &=H_{x} \left(s\right)-\sum _{i,j=1}^{x}\frac{1}{\left(i\cdot j\right)^{s} }  M_{\frac{x^{2} }{i\cdot j} } \left(s\right)\, . 
\end{aligned} 
\end{equation}
Combining \eqref{Eq_3} and \eqref{Eq_4}, we get the identity

\begin{equation} \label{Eq_5}
H_{x^{2} } \left(s\right)=2H_{x} \left(s\right)-\sum _{i,j=1}^{x}\frac{1}{\left(i\cdot j\right)^{s} }  M_{\frac{x^{2} }{i\cdot j} } \left(s\right) .
\end{equation}

Consider two most interesting cases following from \eqref{Eq_5} . At \textit{s} = 0 $M_{x} \left(0\right)=\sum _{k=1}^{x}\mu _{k} \equiv  M_{x} $ is Mertens function and we have
\begin{equation}\label{Eq_6}
\left[ x^{2}\right] =2\left[x\right]-\sum _{i,j=1}^{x}M_{\frac{x^{2} }{i\cdot j} } \,,
\end{equation}
while at \textit{s} = 1 $M_{x} \left(1\right)=\sum _{k=1}^{x}\frac{\mu _{k} }{k} \equiv  \, m_{x} $, the formula \eqref{Eq_5} is an expression for harmonic number at $x^{2} $
\begin{equation}\label{Eq_7}
H_{x^{2} } =2H_{x} -\sum _{i,j=1}^{x}\frac{1}{i\cdot j}  \, m_{\frac{x^{2} }{i\cdot j} } .   
\end{equation}
Applying now the asymptotic formula for harmonic number \cite{Harmonic}, we get
\begin{equation}\label{Eq_8}
\log \left[x^{2} \right]+\gamma +O\left(\frac{1}{x^{2} } \right)=2\left(\log \left[x\right]+\gamma +O\left(\frac{1}{x} \right)\right)-\sum _{i,j=1}^{x}\frac{1}{i\cdot j}  \, m_{\frac{x^{2} }{i\cdot j} } ,
\end{equation}
where  $\gamma $= 0.5772156... is Euler's constant \cite{Constant}. Hence immediately follows the new series for constant $\gamma $
\begin{equation}\label{Eq_9}
\sum _{i,j=1}^{x}\frac{1}{i\cdot j}  \, m_{\frac{x^{2} }{i\cdot j} } =\gamma +O\left(\frac{1}{x} \right)
\end{equation}
or
\begin{equation}\label{Eq_10}
\gamma =\mathop{\lim }\limits_{x\to \infty } \sum _{i,j=1}^{x}\frac{1}{i\cdot j}  \, m_{\frac{x^{2} }{i\cdot j} } .
\end{equation}

By analogy with \eqref{Eq_1}-\eqref{Eq_7}, we can obtain the similar set of equations for generalized harmonic numbers in power \textit{s }$H_{x} \left(s\right)$. In particular, the corresponding counterparts \eqref{Eq_11}-\eqref{Eq_18}, expressed in terms of generalized harmonic numbers, can be found:
\begin{equation}\label{Eq_11}
\sum _{i=1}^{x}\frac{\mu _{i} }{i^{s} }  H_{\frac{x}{i} } \left(s\right)=1 ,
\end{equation}
\begin{equation}\label{Eq_12}
\sum _{i=1}^{x}\frac{\mu _{i} }{i^{s} }  H_{\frac{\left(x^{2} /j\right)}{i} } \left(s\right)=1,\qquad \qquad j\ge x ,
\end{equation}
Multiplying both parts of \eqref{Eq_12} by $\frac{\mu _{j} }{j^{s} } $ and taking sum over index \textit{j} from $x+1$ up to $x^{2} $ yields
\begin{equation}\label{Eq_13}
\sum _{j=x+1}^{x^{2} }\frac{\mu _{j} }{j^{s} } \sum _{i=1}^{x}\frac{\mu _{i} }{i^{s} }   H_{\frac{\left(x^{2} /j\right)}{i} } \left(s\right)=M_{x^{2} } \left(s\right)-M_{x} \left(s\right) ,
\end{equation}

\begin{equation}\label{Eq_14}
\begin{aligned} 
{\sum _{i=1}^{x}\frac{\mu _{i} }{i^{s} }  \sum _{j=x+1}^{x^{2} }\frac{\mu _{j} }{j^{s} }  H_{\frac{\left(x^{2} /i\right)}{j} } \left(s\right)} & {=} & {\sum _{i=1}^{x}\frac{\mu _{i} }{i^{s} }  \left(1-\sum _{j=1}^{x}\frac{\mu _{j} }{j^{s} }  H_{\frac{\left(x^{2} /i\right)}{j} } \left(s\right)\right)} \\ {\qquad \qquad \qquad } & {=} & {M_{x} \left(s\right)-\sum _{i,j=1}^{x}\frac{\mu _{i} \mu _{j} }{\left(i\cdot j\right)^{s} }  H_{\frac{x^{2} }{i\cdot j} } \left(s\right)\, \, .\, } 
\end{aligned} 
\end{equation} 
By analogy with \eqref{Eq_5} we get the identity
\begin{equation}\label{Eq_15}
M_{x^{2} } \left(s\right)=2M_{x} \left(s\right)-\sum _{i,j=1}^{x}\frac{\mu _{i} \mu _{j} }{\left(i\cdot j\right)^{s} }  H_{\frac{x^{2} }{i\cdot j} } \left(s\right)  .
\end{equation}
Consider again two most interesting cases. At \textit{s} = 1 \eqref{Eq_15} gives
\begin{equation}\label{Eq_16}
m_{x^{2} } =2m_{x} -\sum _{i,j=1}^{x}\frac{\mu _{i} \mu _{j} }{i\cdot j}  H_{\frac{x^{2} }{i\cdot j} } ,
\end{equation}
while at \textit{s} = 0 we have
\begin{equation}\label{Eq_17}
M_{x^{2} } =2M_{x} -\sum _{i,j=1}^{x}\mu _{i} \mu _{j}  \left[\frac{x^{2} }{i\cdot j} \right] \end{equation}
or
\begin{equation}\label{Eq_18}
M_{x^{2} } =2M_{x} -x^{2} \cdot m_{x}^{2} +\sum _{i,j=1}^{x}\mu _{i} \mu _{j}  \left\{\frac{x^{2} }{i\cdot j} \right\} .
\end{equation}

Let us consider identity \eqref{Eq_17} in more detail. Assume the Riemann Hypothesis. In this case for all $\varepsilon >0 \ M_{x} =O\left(x^{\frac{1}{2} +\varepsilon } \right)$, which is equivalent to $\left|M_{x} \right|\le C\sqrt{x} \left(\log x\right)^{n} $ for the given positive constants \textit{C} and \textit{n}. For $x^{2} $ this can be rewritten as $\left|M_{x^{2} } \right|\le C\cdot 2^{n} \cdot x\left(\log x\right)^{n} $. Obviously, if $M_{x} =o\left(\sqrt{x} \log x\right)$ then some positive function $f\left(x\right)$ satisfying the conditions $f\left(x\right)=o\left(\log x\right)$ and $M_{x} =O\left(\sqrt{x} f\left(x\right)\right)$ can be used instead of $\log x$.
\\
\\
\indent Let us formulate \textit{the Induction Procedure} for the reverse statement.
\\
\textbf{\textit{\underbar{Induction Procedure} }}
\\
\indent Assume that for some real bound \textbf{$x_{0} >e$}, variables \textit{x} and \textit{y}, and for the given positive constants $C$ and $n$
\begin{equation}\label{Eq_19}
\sup \left|M_{y} \right| =\sup \left|\sum _{i}^{y}\mu _{i}  \right| \le C\sqrt{x} \left(\log x\right)^{n} ,    {\qquad}   e\le y\le x<x_{0} .
\end{equation}

In first step the verification of this assumption can be done by direct calculation of 
$\left|M_{k} \right|$ for all $k<x_{0} $, where $k$ is the natural number.

If applying the identity \eqref{Eq_17} we can prove that always (regardless of \textbf{$x_{0} $} value)
\begin{equation}\label{Eq_20}
\begin{aligned}
\sup \left|M_{y^{2} } \right| &\sim \sup \left|\sum _{i,j=1}^{y}\mu _{i} \mu _{j}  \left[\frac{y^{2} }{i\cdot j} \right]\right|\\&\le 2^{n} \cdot \sqrt{x} \cdot \sup \left|M_{y} \right| , {\qquad \qquad \qquad } e\le y\le x<x_{0}  ,
\end{aligned}
\end{equation}
then we confirm that
$$
\sup \left|M_{y} \right| \sim \sup \left|\sum _{i,j=1}^{\sqrt{y} }\mu _{i} \mu _{j}  \left[\frac{y}{i\cdot j} \right]\right|\\ \le C\sqrt{x} \left(\log x\right)^{n} , {\qquad } e\le y\le x<x_{0}^{2}  ,
$$
and the statement \eqref{Eq_19} is extended now up to $x_{0}^{2} $. 
\\
\textbf{\textit{\underbar{End of Induction Procedure}}}
\\
\\
\textit{The Induction Procedure} can be applied over and over again for further validation of \eqref{Eq_19}. Hence the Riemann Hypothesis is justified.
\\
\\
Thus, to prove the Riemann Hypothesis it is enough to prove that: 

\textbf{if for some real bound $x_{0}$, variables \textit{x} and \textit{y}, and for the \indent given positive constants }$C$ \textbf{and} $n$\textbf{}
\begin{equation}\label{Eq_21}
\sup \left|\sum _{i}^{y}\mu _{i}  \right|\le C\sqrt{x} \left(\log x\right)^{n},  {\qquad }    e\le y\le x<x_{0} ,
\end{equation}

\textbf{then always (independently of $x_{0} $ value) follows}

\begin{equation}\label{Eq_22}
\begin{aligned}
& \sup \left|\sum _{i,j=1}^{y}\mu _{i} \mu _{j}  \left[\frac{y^{2} }{i\cdot j} \right]\right|= \sup \left|y^{2} m_{y}^{2} -\sum _{i,j=1}^{y}\mu _{i} \mu _{j}  \left\{\frac{y^{2} }{i\cdot j} \right\}\right| \\ &\le 2^{n} \cdot \sqrt{x} \cdot \sup \left|\sum _{i}^{y}\mu _{i}  \right|, {\qquad } e\le y\le x<x_{0} .
\end{aligned}
\end{equation}
\textbf{}

\bigskip
\bigskip
\bigskip


\begin{thebibliography}{8}


\bibitem{RAbrarov}
R. M. Abrarov, S. M. Abrarov, \textit{On the properties of generalized harmonic and oscillatory numbers. Simple proof of the Prime Number Theorem}, \indent\verb+http://arxiv.org/PS_cache/arxiv/pdf/0709/0709.3145v2.pdf+

\bibitem{Harmonic}
\indent\verb+http://mathworld.wolfram.com/HarmonicNumber.html+


\bibitem{Constant}
\indent\verb+http://mathworld.wolfram.com/Euler-MascheroniConstant.html+
\\
\end{thebibliography}
\end{document}